\documentclass[12pt]{article}
\usepackage{mathrsfs}

\usepackage{amsfonts}
\usepackage{}

\usepackage{amssymb}
\usepackage{amsmath}
\usepackage{graphics}
\usepackage{epsfig}
\newtheorem{claim}{\bf \t}[part]

\newtheorem{Corollary}{Corollary}[part]
\newtheorem{Definition}{Definition}[part]

\newtheorem{Lemma}{Lemma}[part]
\newtheorem{Proposition}{Proposition}[part]
\newtheorem{Remark}{Remark}[part]
\newtheorem{Theorem}{Theorem}[part]

\numberwithin{Assumption}{section} \numberwithin{Corollary}{section}
\numberwithin{Definition}{section} \numberwithin{equation}{section}
\numberwithin{Example}{section} \numberwithin{Lemma}{section}
\numberwithin{Proposition}{section} \numberwithin{Remark}{section}
\numberwithin{Theorem}{section}

\def\t{\theta}

\def\f{\frac}

\def\text#1{{\rm #1}}

\begin{document}
\date{}

\title{\Large \bf  Blow-up of Smooth Solutions to the Euler-Poisson Equations}

\author{\small \textbf{Yuexun Wang}\thanks{ E-mail:
yx-wang13@mails.tsinghua.edu.cn}}
 \maketitle \small
Mathematical Sciences center, Tsinghua University, Beijing
100084, P. R. China

\small \maketitle { \bf Abstract:}  In this paper, the finite time blow-up of
 smooth solutions to the Cauchy problem for full
Euler-Poisson equations and  isentropic Euler-Poisson equations with repulsive forces or attractive  forces in high dimensions $(n\geq3)$ is proved for a large class of initial data. It is  not
required that the initial data has compact support or contains vacuum
in any finite regions.

\section{Introduction } \setcounter{equation}{0}
\setcounter{Assumption}{0} \setcounter{Theorem}{0}
\setcounter{Proposition}{0} \setcounter{Corollary}{0}
\setcounter{Lemma}{0}

    \qquad The purpose of this paper is to prove the finite time blow-up of the system of compressible Euler-Poisson equations in multi-dimensions for both the gravitational forces modeling the evlotion of gaseous stars and  the repulsive forces modeling the semiconductor devices.  For the compressible Euler equations in three dimensions, the finite time blow-up of smooth solutions was proved by Sideris in \cite{Si}, with the assumption that the difference of the initial data and a constant state has a compact support. For the compressible Navier-Stokes equations with vacuuum, Xin in \cite{X} proved the finite blow-up of smooth solutions for the initial data with compact suppport. The arguments in \cite{Si} and \cite{X} both rely on the finite propagation speed of the supports of solutions. However, for the coppressible Euler-Poission equations with the  gravitational forces or repulsive forces, the gravitational or electric potential are non-local, satisfying elliptic equations, and thus the compressible Euler-Poisson equations do not possess the property of the finite propagation speed. This makes the study of the finite time blow up of smooth solutions challenging and interesting. For the spherically symmetric solutions, it was proved by Makino and Pertheme in \cite{MP} that the smooth solutions blow up in a finite time for the isentropic Euler-Poisson equations with grivatational forces,  if the initial data is symmetric  and has a  compact support. The argument in \cite{MP} uses the sign of the gravitational force in the
spherical symmetry case. In this paper, we prove the blow-up of solutions to the general multi-dimensional  compressible Euler-Poisson equations.
The system of full Euler-Poisson equations in $\mathbb{R}^n(n\geq3)$ is the following:
\begin{eqnarray}\label{(EP-1)}
\left\{ \begin{array}{ll}
\partial_t\rho+\textrm{div}(\rho
u)=0,\\
\partial_t(\rho u)+\textrm{div}(\rho u \otimes u)+\nabla p=\delta\rho \nabla \Phi,\\
\partial_t({\f12}\rho |u|^2+\rho e )+\textrm{div}(({\f12}\rho |u|^2+\rho e+p)u)=0,\\
\Delta \Phi= n(n-2)\omega_ng\rho.
\end{array}
\right.
\end{eqnarray}
Here $(x,t)\in\mathbb{R}^n\times\mathbb{R}_+ $, $\rho= \rho
(x,t), u = (u_1, u_2,\cdots ,u_n), p$, $e$ and $\Phi $   denote the
density, velocity, pressure, internal
energy and  the potential for the self-gravitational force, respectively. The constants $g$
and $\omega_n$ are gravitational constant and the measure of the unite ball in $\mathbb{R}^n$, respectively.

 The polytropic gas satisfies the following state  equations:
\begin{equation}\label{(1)}
 p=(\gamma-1)\rho e, \ {\rm or}\ ,    p =Aexp({\f{s}{c_\nu}})\rho^\gamma,
\end{equation}
where  $R > 0$ is the gas constant, $A>0$ is an absolute constant,
$\gamma> 1$  is the specific heat ratio, $c_\nu= {\f{R}{\gamma-1}}$
and $s$ is the entropy.

The initial data to  the equations \eqref{(EP-1)} are imposed as
 \begin{equation}\label{(2)}
  (\rho,u,s)(x,t)|_{t=0}=(\rho_0(x),u_0(x),s_0(x)).
 \end{equation}

If the entropy is a constant, the full Euler-Poisson equations
 reduce to isentropic ones which read as
\begin{eqnarray}\label{(EP-2)}
\left\{ \begin{array}{ll}
\partial_t\rho+\textrm{div}(\rho
u)=0,\\
\partial_t(\rho u)+\textrm{div}(\rho u \otimes u)+\nabla p=\delta\rho \nabla \Phi,\\
\Delta \Phi= n(n-2)\omega_ng\rho.
\end{array}
\right.
\end{eqnarray}

The state  equation of the isentropic process becomes
\begin{equation}\label{(3)}
  p=A\rho^\gamma.
 \end{equation}

The initial data to the equations \eqref{(EP-2)} are imposed as

\begin{equation}\label{(4)}
  (\rho,u)(x,t)|_{t=0}=(\rho_0(x),u_0(x)).
 \end{equation}

The  equation $(1.1)_4$ or the  equation $(1.4)_3$
can be solved as
\begin{equation}\label{(5)}
  \Phi(x,t)=\int_{\mathbb{R}^n}G(x-y)\rho(y,t)dy,
 \end{equation}
where G is the Green's function for the Poisson equation defined by
\begin{eqnarray}\label{(6)}
G(x)={\f{-1}{|x|^{n-2}}}, n\geq3.
\end{eqnarray}

When $\delta=-1$, the system is self-gravitating. The equations (2) are the Newtonian descriptions of
gaseous stars or a galaxy in astrophysics \cite{BT} and \cite{Cha}. When $\delta=1$, the system is the compressible
Euler-Poisson equations with repulsive forces. It can be used as a semiconductor model \cite{Che}.

The results for existence theories and stabilities can be founded in  \cite{BE}, \cite{CT}, \cite{DLYY}, \cite{Ga}, \cite{Gu},
\cite{J}, \cite{LL}, \cite{Li}, \cite{Lin}, \cite{LT1}, \cite{LT2}, \cite{LT3}, \cite{LT4}, \cite{LS1}, \cite{LS2}, \cite{LS3}, \cite{M2}, \cite{M3}, \cite{TW}, \cite{T},  \cite{Wei}
 and  \cite{WTH}. On the other hand,
for the construction of the analytical solutions for the system with pressure,
the interested readers may see \cite{Cha}, \cite{DXY}, \cite{Go}, \cite{M2}, \cite{Y1} and \cite{Y2}.

 The study of singularity formation in the model with diffusion and relaxation in one-dimensional case can be found in \cite{WC} by using the characteristic method.
Local conditions for the finite-time loss of smoothness in the one-dimensional case
with and without pressure were given in \cite{E} and \cite{ELT}.  The non-blowup phenomenon was
discussed  in \cite{DXY}. Recently,    the finite time blowup results
for the pressureless case were given in \cite{CT},  \cite{LL} and \cite{Y2}. The author  studies finite time blowup  if the  total energy is non-positive in \cite{Ye}.

For  the isentropic Euler-Poisson equations with repulsive forces.  The non-existence of global gentle solutions has also been proved in \cite{P}  if the initial data  is symmetric with compact support.

Our aim in this paper is to show that finite time blowup is generic for  the Cauchy problem to  full
Euler-Poisson equations and  isentropic Euler-Poisson equations with repulsive forces or attractive forces in high dimensions $(n\geq3)$ under suitable conditions for a large class of initial data. It is  not
required that the initial data has compact support or contains vacuum
in any finite regions.
The method is to compare the upper  and lower bound estimates on the internal energy. To obtain our main results,
some physical quantities such as mass, momentum, momentum of
inertia, internal energy, potential energy, total energy and some
combined functionals of these quantities are introduced. The basic
properties and specific  relationships between them are helpful and
crucial to prove the main results.  More precisely, for the full
Euler-Poisson equations, the upper and lower decay rates
of the internal energy will be calculated in a precise way. And for
the isentropic Euler-Poisson equations, the upper and
lower decay rates of the potential energy will be presented
accordingly. Then the main results will be proved by comparing the
coefficients of the upper and lower decay rates.

In the following, we denote the full Euler-Poisson equations \eqref{(EP-1)} and \eqref{(1)} as (EP),
the isentropic Euler-Poisson equations \eqref{(EP-2)} and \eqref{(3)} as (IEP) for short.
And for convenience, we take $g=1$ and $A=1$.

The paper is organized as follows. In section 2, we will introduce some physical quantities and functionals and present our main results. In Section 3, we will give some basic properties of the physical quantities and functionals and the proof of the main results.

\section{Main Results }

The following  physical quantities will be  used in this paper.
\begin{equation}\label{(Q1)}
\left\{
\begin{array}{lll}
&M(t)=\int _{\mathbb{R}^n}\rho  dx \ \ &(mass),\\
&\mathbb{P}(t)=\int_{\mathbb{R}^n} \rho u dx\ \ &(momentum),\\
&F(t)=\int_{\mathbb{R}^n} \rho u x dx &(momentum\   weight), \\
&G(t)={\f12}\int_{\mathbb{R}^n} \rho |x|^2 dx \ \ &(momentum \ of \ inertia),\\
&E_{\delta}(t)={\f12}\int _{\mathbb{R}^n} \rho |u|^2 dx+\int_{\mathbb{R}^n}
\rho e dx-{\f{\delta}{2}} \int_{\mathbb{R}^n}
\rho \Phi dx\\
& \quad\quad\triangleq
 E_k(t)+E_i(t)+E_p(t)\ \ &(total \ energy),\\
&IE_{\delta}(t) ={\f12}\int_{\mathbb{R}^n} \rho |u|^2 dx+{\f{1}{\gamma
-1}}\int_{\mathbb{R}^n} P dx-{\f{\delta}{2}} \int_{\mathbb{R}^n}
\rho \Phi dx\\
&\quad\quad\triangleq E_k(t)+I(t)+E_p(t) \ \ &(energy).
\end{array}
\right.
\end{equation}

The basic properties and relationships among these quantities will
be discussed in Section 3. To prove our main results, we
introduce the following functionals:
\begin{equation}\label{(Q2)}
\left\{
\begin{array}{ll}
&H_\delta(t)=2 E_k(t)+n(\gamma -1)E_i(t)-{\f{\delta(n-2)}{2}} \int_{\mathbb{R}^n}
\rho \Phi dx,\\
&IH_\delta(t)=2 E_k(t)+n(\gamma
-1)I(t)-{\f{\delta(n-2)}{2}} \int_{\mathbb{R}^n}
\rho \Phi dx,\\
&J_\delta(t)= G(t)-(t+1)F(t)+(t+1)^2 E_\delta(t),\\
& IJ_\delta(t)= G(t)-(t+1)F(t)+(t+1)^2 IE_\delta(t).
\end{array}
 \right.
\end{equation}

 For any $T>0$, we require $\rho, u$ and $ s$ to satisfy
\begin{eqnarray}\label{(1001)}
\rho|x|^2, p, \rho |u|^2, \rho \Phi \in L^\infty((0,T);L^1(\mathbb{R}^n)), \nabla u\in L^2((0,T);L^2(\mathbb{R}^n)).
\end{eqnarray}

We always assume that $M(0), \mathbb{P}(0), F(0), G(0),  E_\delta(0),
IE_\delta(0)$ are finite and   $M(0)>0$. Since $F(t)^2\leq4G(t)E_k(t)$ (see Lemma \ref{(100000)}) and
$E_i(t)\ge 0, I(t)\ge 0$, it follows that $J_\delta(0)\geq0$ and $IJ_\delta(0)\geq0$,
respectively.

It should be remarked that the conditions
\eqref{(1001)}  guarantee that the integration by parts
in our calculations make sense. Actually, \eqref{(1001)} is equivalent to $G(t), E_{\delta}(t)$ and $IE_{\delta}$ are finite on $[0,T]$.
 In particular, for the classical solution to  (EP)
satisfying \eqref{(1001)}, the conservations of the
mass $M(t)$,  momentum $\mathbb{P}(t)$ and  energy $E(t)$ hold true.
Therefore we naturally introduce the following two definitions.

\begin{Definition}\label{(1000)}
For the Cauchy problem to (EP), we call $ (\rho,u,s)\in X(T)$ if $ (\rho,u,s)$ is a
classical solution in $[0,T]$  and satisfies \eqref{(1001)}.
\end{Definition}

\begin{Definition}\label{(1002)}
For the Cauchy problem to (IEP), we call
$ (\rho,u)\in Y(T)$ if $ (\rho,u)$ is a classical solution in
$[0,T]$ and  satisfies \eqref{(1001)}.
\end{Definition}

In order to present precisely our results  freely, we  need to introduce some symbols to denote the numbers used later:

$$C_{HLP}---Hardy-Littlewood-Paley \ {\rm constant}\ ,$$
$$C_{HLS}---Hardy-Littlewood-Sobolev \ {\rm constant}\ ,$$
$$ s_1=\min_x s_0(x),C_0= IE_\delta(0), C_1=\min_{p,q} C_{HLS}M(0)^{2-\theta}, $$
\begin{eqnarray}\nonumber
C_2=
\left\{
\begin{array}{ll}
&2n^{2+{\f{n-2}{2+n(\gamma-1)}}}(n-2)^{1+{\f{n-2}{2+n(\gamma-1)}}}\omega_n^{2+{\f{n-2}{2+n(\gamma-1)}}}\times\\
&M^{2+{\f{(n-2)(2-\gamma)}{2+n(\gamma-1)}}}C_{HLP}^{\f{(n-2)(2-\gamma)}{2+n(\gamma-1)}},\quad \quad\quad \quad\quad \quad\quad 2(1-{\f{1}{n }})<\gamma<2,\\
& 2M(\gamma-2)+2M^2n^{1+{\f{n}{2}}}(n-2)^{\f{n}{2}}\omega_n^{1+{\f{n}{2}}}(\gamma-1)^{1-{\f{n}{2}}},\gamma\geq2,
\end{array}
 \right.
\end{eqnarray}
$$C_3=\max\{6-n,n(2\gamma-3)+2\}C_0+\max\{4-n,n(\gamma-2)+2\}C_2,$$
$$C_4=\max\{2,n(\gamma-1)\}C_0, C_5=\max\{2,n(\gamma-1)\}C_0,$$
$$C_6=\min\{4-n,n(\gamma-2)+2\}(E_k(0)+E_i(0))+(n-2)E_\delta(0),$$
$$C_7=\max\{4-n,n(\gamma-2)+2\}(E_k(0)+E_i(0))+(n-2)E_\delta(0),$$
$$C_8=2({\f{\pi^{\f n 2}}{\Gamma ({\f n 2}+1)}})^{\f{2(\gamma-1)}{(n+2)\gamma-n}},C_{9}=({\f{\Gamma ({\f n 2}+1)}{(\pi)^{\f n
2}}})^{\gamma-1 } {\f{\exp ({\f {s_1}{c_\nu}})M(0)^{\f
{(n+2)\gamma-n}{2}}}{2^{\f{(n+2)\gamma-n}{2}}(\gamma-1)}},$$
$$C_{10}=({\f{\Gamma ({\f n 2}+1)}{(\pi)^{\f n 2}}})^{\gamma-1 }
{\f{M(0)^{\f{(n+2)\gamma-n}{2}}}{2^{\f{(n+2)\gamma-n}{2}}(\gamma-1)}}, C_{11}=IJ_\delta(0).$$

Our main results  are presented as follows.
 \begin{Theorem}\label{(10000)}
  For (IEP) with $\delta=-1$. If  the initial values satisfy one of the following conditions:\\
  (i) $\gamma >2(1-{\f{1}{n }})$ and
  \begin{equation}\label{(7)}
   C_3<0,
\end{equation}
  (ii) $\gamma >2(1-{\f{1}{n }})$ and
  \begin{equation}\label{(8)}
   C_3=0, F(0)<0,
\end{equation}
  (iii) ${\f43}<\gamma\leq{\f53}$, $n=3$ and
  \begin{equation}\label{(9)}
   (C_3\triangleq)3C_0+C_2>0, C_{10}>C_{11}(3C_0+C_2)^{\f{3(\gamma-1)}{2}},
\end{equation}
 then there   is no solution in $Y(\infty)$ to the
Cauchy problem of (IEP).  And, in the case of (iii), there exists a time $T_1^{*}>0$
such that there is no solution in $Y(T_1^{*})$.
\end{Theorem}
\begin{Theorem}\label{(10001)}
  For (IEP) with $\delta=1$. Assume $1<\gamma \leq1+{\f{2}{n }}$ and $n\geq4$. If  the initial values satisfy  the following condition:
  \begin{equation}\label{(10)}
  C_{10}>2^{\f{n(\gamma-1)}{2}}C_{11},
\end{equation}
 then there exists a time $T_2^{*}>0$
such that there is no solution in $Y(T_2^{*})$ to the
Cauchy problem of (IEP).
\end{Theorem}

\begin{Theorem}\label{(10002)}
 For (EP) with $\delta=-1$. If  the initial values satisfy one of the following conditions:\\
 (i) \begin{equation}\label{(11)}
  C_7<0,
\end{equation}
 (ii)\begin{equation}\label{(12)}
  C_7=0, F(0)<0,
\end{equation}
then   there is no solution in
$X(\infty)$ to the
Cauchy problem of (EP).
\end{Theorem}

\begin{Remark}
Global
existence for the
Cauchy problem to (IEP) with $\delta=1$ was obtained by Guo in \cite{Gu}, assuming the flow is irrotational and the data is
in the small $H^2$-neighborhood of a constant state.
\end{Remark}

\begin{Remark}
 In the system of  full Euler-Poisson equations, the energy equation can be replaced by the entropy equation. From
 the entropy equation, one can deduce that the entropy $s(x,t)$ is increasing in time. Therefore we conclude that
 $s(x,t)\geq s_0(x)$.
\end{Remark}

\begin{Remark}
In \eqref{(7)}-\eqref{(12)}, we do not  require that the
initial data has compact support or contains  vacuum in any finite
region.
\end{Remark}

\begin{Remark}
 The time $T_1^{*}$ and $T_2^{*}$ can be computed precisely (see the proof of Theorem \ref{(10000)}). In other words, we can
 find out the "last" blow-up time.
\end{Remark}

\section{ The Proof of Theorem \ref{(10000)}-\ref{(10002)}}

Let $(\rho,u,s)\in X(T) $ be a classical solution to the Cauchy
problem to (EP). And $(\rho,u)\in Y(T) $ is a classical
solution to the Cauchy problem  to (IEP). In this
subsection, we will first present some basic relationships among
the quantities defined in Section 2. Then we will give the proof of
Theorem \ref{(1000)}-Theorem \ref{(1001)}.

We first begin with the  relationships among physical quantities and functionals listed in Section 2.

\begin{Lemma}\label{(100000)}
For (EP) and (IEP), we have
\begin{eqnarray}\label{(13)}
{\f {d}{dt}}M(t)=0,  {\f {d}{dt}}\mathbb{P}(t)=\delta \int_{\mathbb{R}^n}
\rho \nabla\Phi dx,  {\f {d}{dt}}G(t)=F(t).
\end{eqnarray}
For (EP), we have
\begin{eqnarray}\label{(14)}
{\f {d}{dt}}E_\delta(t)=-{\f{\delta}{2}} {\f {d}{dt}}\int_{\mathbb{R}^n}
\rho \Phi dx,  {\f {d}{dt}}F(t)=H_\delta(t).
\end{eqnarray}
For (IEP), we have
\begin{eqnarray}\label{(15)}
{\f {d}{dt}}IE_\delta(t)=0,{\f {d}{dt}}F(t)=IH_\delta(t).
\end{eqnarray}
\end{Lemma}
 \textbf{{\em Proof.}} Using (EP) and (IEP), applying integration by
 parts,  one can verify \eqref{(13)}-\eqref{(15)}.

 By $H\ddot{o}lder's$ inequality, we  arrive at the following:

\begin{Lemma}\label{(100001)}
 For (EP) and (IEP), we have
\begin{equation}\label{(16)}
F(t)^2\leq4G(t)E_k(t).
\end{equation}
\end{Lemma}

We next quote two famous inequalities. One is  Hardy-Littlewood-Sobolev inequality, which can be presented as follows:

\begin{Lemma}\label{(100002)}
For all $f\in L^p(\mathbb{R}^n)$, $g\in L^q(\mathbb{R}^n)$, $1<p,q<\infty$, $0<\lambda<n$ and ${\f{1}{p}}+{\f{1}{q}}+{\f{\lambda}{n}}=2$, it holds
 \begin{equation}\label{(17)}
|\int_{\mathbb{R}^n}\int_{\mathbb{R}^n}f(x)|x-y|^{-\lambda}g(y)dxdy|\leq C_{HLS}||f||_{L^p(\mathbb{R}^n)}
||g||_{L^q(\mathbb{R}^n)},
\end{equation}
where $C_{HLS}={\f{1}{pq}}{\f{n}{n-\lambda}}({\f{\omega_{n-1}}{n}})^{\lambda/n}
(({\f{\lambda/n}{1-1/p}})^{\lambda/n}+({\f{\lambda/n}{1-1/q}})^{\lambda/n})$.
\end{Lemma}

The other is  Hardy-Littlewood-Paley inequality, which  reads as:
\begin{Lemma}\label{(100004)}
 If $1<p\leq2$ and $f\in L^p(\mathbb{R}^n)$, then there exists a positive constant $C_{HLP}$ such that
 \begin{equation}\label{(22)}
(\int_{\mathbb{R}^n}|\mathscr{F}f(\xi)|^p|\xi|^{n(p-2)}d\xi)^{\f{1}{p}}\leq C_{HLP}||f||_{L^p(\mathbb{R}^n)},
\end{equation}
where $\mathscr{F}f(\xi)=\int_{\mathbb{R}^n}f(x)e^{-2\pi i x\xi}dx$.
\end{Lemma}

From the two famous inequalities above, we can obtain the following two  estimates, which are very useful to control the lower
bound and upper bound of internal energy.

\begin{Lemma}\label{(100003)}
If $\gamma>{\f{2n}{n+2}}$, then
 \begin{equation}\label{(18)}
-\int_{\mathbb{R}^n}\rho \Phi dx=\int_{\mathbb{R}^n}\int_{\mathbb{R}^n}\rho(x)|x-y|^{2-n}\rho(y)dxdx\leq C_1
||\rho||_{L^\gamma(\mathbb{R}^n)}^\theta,
\end{equation}
where $\theta={\f{(n-2)\gamma}{n(\gamma-1)}}\in (0,2)$, $C_1=\min_{p,q} C_{HLS}M(0)^{2-\theta}$ and $\min_{p,q} C_{HLS}$ is the minimum of $C_{HLS}$ for choosing some $p$ and $q$.
\end{Lemma}
\textbf{{\em Proof.}}  Young's inequality yields
\begin{eqnarray}\label{(19)}
||\rho||_{L^p(\mathbb{R}^n)}\leq ||\rho||_{L^1(\mathbb{R}^n)}^{1-\alpha}||\rho||_{L^\gamma(\mathbb{R}^n)}^{\alpha},
\end{eqnarray}
where ${\f{1}{p}}=1-\alpha+{\f{\alpha}{\gamma}}$.
And
\begin{eqnarray}\label{(20)}
||\rho||_{L^q(\mathbb{R}^n)}\leq ||\rho||_{L^1(\mathbb{R}^n)}^{1-\beta}||\rho||_{L^\gamma(\mathbb{R}^n)}^{\beta}.
\end{eqnarray}
where ${\f{1}{q}}=1-\beta+{\f{\beta}{\gamma}}$.
Therefore, we have
\begin{eqnarray}\label{(21)}
||\rho||_{L^p(\mathbb{R}^n)}||\rho||_{L^q(\mathbb{R}^n)}\leq ||\rho||_{L^1(\mathbb{R}^n)}^{2-(\alpha+\beta)}||\rho||_{L^\gamma(\mathbb{R}^n)}^{\alpha+\beta}
=C_{HLS}M(0)^{2-(\alpha+\beta)}||\rho||_{L^\gamma(\mathbb{R}^n)}^{\alpha+\beta}.
\end{eqnarray}
Now if we define  $\theta\triangleq\alpha+\beta={\f{(n-2)\gamma}{n(\gamma-1)}}$, then $\theta\in (0,2)$ since ${\f{1}{p}}+{\f{1}{q}}+{\f{\lambda}{n}}=2$ and $\gamma>{\f{2n}{n+2}}$.
This, together with Lemma \ref{(100002)}, implies \eqref{(18)}.

Employing the trick in \cite{DLYY} and \cite{DXY}, we arrive at the following estimates.

\begin{Lemma}\label{(100005)}
 For (EP) and (IEP), if $\gamma>2(1-{\f{1}{n }})$, then for any $\varepsilon>0$, there exists a positive constant $C(\varepsilon,n,\omega_n,M)$ such that
 \begin{equation}\label{(23)}
 -\int_{\mathbb{R}^n}\rho \Phi dx\leq \varepsilon E_{i}(t)+C(\varepsilon,n,\omega_n,M),
\end{equation}
and
\begin{equation}\label{(24)}
 -\int_{\mathbb{R}^n}\rho \Phi dx\leq \varepsilon I(t)+C(\varepsilon,n,\omega_n,M).
\end{equation}
\end{Lemma}
\textbf{{\em Proof.}} We only treat with the estimate \eqref{(24)} corresponding to (IEP), since  \eqref{(23)} is similar.
We first consider the case of  $2(1-{\f{1}{n }})<\gamma<2$.
 Taking Fourier transformation on  $(1.4)_3$ leads to
 \begin{equation}\label{(24.5)}
 \mathscr{F}\Phi(\xi)=-{\f{n(n-2)\omega_n}{4\pi^2}}\mathscr{F}\rho(\xi)|\xi|^{-2}.
 \end{equation}
 Then for any $r>0$,  by Plancherel's theorem, we get
\begin{eqnarray}\label{(25)}
&&-\int_{\mathbb{R}^n}\rho \Phi dx
=-{\f{1}{n(n-2)\omega_n}}\int_{\mathbb{R}^n}\Delta\Phi \Phi dx
={\f{1}{n(n-2)\omega_n}}\int_{\mathbb{R}^n}|\nabla \Phi|^2 dx\nonumber\\
&=&{\f{1}{n(n-2)\omega_n}}\int_{\mathbb{R}^n}|\mathscr{F}(\nabla \Phi)|^2 d\xi
=n(n-2)\omega_n\int_{\mathbb{R}^n}|\mathscr{F}\rho(\xi)|^2|\xi|^{-2} d\xi\nonumber\\
&=&n(n-2)\omega_n(\int_{|\xi|>r}|\mathscr{F}\rho(\xi)|^2|\xi|^{-2} d\xi+\int_{|\xi|\leq r}|\mathscr{F}\rho(\xi)|^2|\xi|^{-2} d\xi)\nonumber\\
&\leq&n(n-2)\omega_n(M^{2-\gamma}r^{-2-n(\gamma-2)}\int_{|\xi|>r}|\mathscr{F}\rho(\xi)|^\gamma|\xi|^{n(\gamma-2)} d\xi+M^2r^{n-3}\int_{|\xi|\leq r}|\xi|^{1-n} d\xi)\nonumber\\
&\leq&n(n-2)\omega_nM^{2-\gamma}r^{-2-n(\gamma-2)}C_{HLP}^\gamma I(t)+n^2(n-2)\omega_n^2M^2r^{n-2}.
\end{eqnarray}
Where we have used Lemma \ref{(100004)} and $||\mathscr{F}\rho||_{L^\infty(\mathbb{R}^n)}\leq||\rho||_{L^1(\mathbb{R}^n)}$.

Next we  consider the other case of $\gamma\geq2$. For any $r>0$, by interpolation inequality, it holds that
\begin{eqnarray}\label{(26)}
&&-\int_{\mathbb{R}^n}\rho \Phi dx
=n(n-2)\omega_n\int_{\mathbb{R}^n}|\mathscr{F}\rho(\xi)|^2|\xi|^{-2} d\xi\nonumber\\
&=&n(n-2)\omega_n(\int_{|\xi|>r}|\mathscr{F}\rho(\xi)|^2|\xi|^{-2} d\xi+\int_{|\xi|\leq r}|\mathscr{F}\rho(\xi)|^2|\xi|^{-2} d\xi)\nonumber\\
&\leq&n(n-2)\omega_nr^{-2}\int_{\mathbb{R}^n}\rho(x)^2 dx+n^2(n-2)\omega_n^2M^2r^{n-2}\nonumber\\
&\leq&n(n-2)\omega_nr^{-2}({\f{\gamma-2}{\gamma-1}}M+{\f{1}{\gamma-1}}\int_{\mathbb{R}^n}\rho(x)^\gamma dx)+n^2(n-2)\omega_n^2M^2r^{n-2}\nonumber\\
&\leq&{\f{n(n-2)\omega_nr^{-2}}{\gamma-1}}I(t)+{\f{n(n-2)\omega_nMr^{-2}(\gamma-2)}{\gamma-1}}
+n^2(n-2)\omega_n^2M^2r^{n-2}.
\end{eqnarray}
Therefore,  by \eqref{(25)} and \eqref{(26)}, for any $\varepsilon>0$, we can choose appropriate $r>0$ such that \eqref{(24)} holds.

\begin{Remark}
Under a stronger condition $\gamma>{\f{n}{2}}$, we can easily prove  a weaker version of Lemma \ref{(100005)}. In fact,   Lemma \ref{(100003)} and Young's inequality
imply \eqref{(23)} and \eqref{(24)} since $\theta<1$ at this moment.
\end{Remark}

\begin{Corollary}\label{(100006)}
 (I) For (IEP) with $\delta=-1$.\\
 (i) If $\gamma>2(1-{\f{1}{n }})$, then
\begin{eqnarray}\label{(27)}
E_k(t)+I(t)\leq2C_0+C_2.
\end{eqnarray}
where
\begin{eqnarray}\label{28}\nonumber
C_2=
\left\{
\begin{array}{ll}
&2n^{2+{\f{n-2}{2+n(\gamma-1)}}}(n-2)^{1+{\f{n-2}{2+n(\gamma-1)}}}\omega_n^{2+{\f{n-2}{2+n(\gamma-1)}}}\times\\
&M^{2+{\f{(n-2)(2-\gamma)}{2+n(\gamma-1)}}}C_{HLP}^{\f{(n-2)(2-\gamma)}{2+n(\gamma-1)}},\quad \quad\quad \quad\quad \quad\quad 2(1-{\f{1}{n }})<\gamma<2,\\
& 2M(\gamma-2)+2M^2n^{1+{\f{n}{2}}}(n-2)^{\f{n}{2}}\omega_n^{1+{\f{n}{2}}}(\gamma-1)^{1-{\f{n}{2}}},\gamma\geq2.
\end{array}
 \right.
\end{eqnarray}
(ii) It holds that
\begin{eqnarray}\label{(29)}
E_k(t)+I(t)\geq C_0.
\end{eqnarray}
(II) For (IEP) with $\delta=1$. We have
\begin{eqnarray}\label{(29.5)}
E_k(t)+I(t)\leq C_0.
\end{eqnarray}
(III) For (EP). It holds that
 \begin{eqnarray}\label{(32)}
E_k(t)+E_i(t)= E_k(0)+E_i(0).
\end{eqnarray}
\end{Corollary}
\textbf{{\em Proof.}}  In view of Lemma \ref{(100000)} and Lemma  \ref{(100005)} with $\varepsilon=1$,
we have
\begin{eqnarray}\label{(30)}
E_k(t)+I(t)
= C_0-{\f12}\int_{\mathbb{R}^n}\rho \Phi dx
\leq C_0+{\f12}I(t)+C(n,\omega_n,M),
\end{eqnarray}
where
\begin{eqnarray}\label{(31)}
C(n,\omega_n,M)=
\left\{
\begin{array}{ll}
&n^{2+{\f{n-2}{2+n(\gamma-1)}}}(n-2)^{1+{\f{n-2}{2+n(\gamma-1)}}}\omega_n^{2+{\f{n-2}{2+n(\gamma-1)}}}\times\\
&M^{2+{\f{(n-2)(2-\gamma)}{2+n(\gamma-1)}}}C_{HLP}^{\f{(n-2)(2-\gamma)}{2+n(\gamma-1)}},\quad \quad\quad \quad\quad \quad\quad 2(1-{\f{1}{n }})<\gamma<2,\\
& M(\gamma-2)+M^2n^{1+{\f{n}{2}}}(n-2)^{\f{n}{2}}\omega_n^{1+{\f{n}{2}}}(\gamma-1)^{1-{\f{n}{2}}},\gamma\geq2.
\end{array}
 \right.
\end{eqnarray}
This implies \eqref{(27)} by taking $C_2=2C(n,\omega_n,M)$. \eqref{(29)} and \eqref{(29.5)} hold obviously since $\int_{\mathbb{R}^n}\rho \Phi dx<0$. \eqref{(32)} follows from Lemma \ref{(100000)}.

We have the following estimates of G(t), which are key parts of obtaining the lower bound of internal energy.

\begin{Proposition}\label{(100008)}
 (I) For (IEP) with $\delta=-1$.\\
 (i) If $\gamma>2(1-{\f{1}{n }})$, then
\begin{eqnarray}\label{(33)}
G(t)\leq C_3t^2+F(0)t+G(0),
\end{eqnarray}
where $C_3=\max\{6-n,n(2\gamma-3)+2\}C_0+\max\{4-n,n(\gamma-2)+2\}C_2$.\\
(ii) It holds that
\begin{eqnarray}\label{(34)}
G(t)\geq C_4t^2+F(0)t+G(0),
\end{eqnarray}
where $C_4=\min\{2,n(\gamma-1)\}C_0$.\\
(II) For (IEP) with $\delta=1$. We have
\begin{eqnarray}\label{(35)}
G(t)\leq C_5t^2+F(0)t+G(0).
\end{eqnarray}
where $C_5=\max\{2,n(\gamma-1)\}C_0$.\\
(III) For (EP). It holds that
\begin{eqnarray}\label{(37)}
 C_6t^2+F(0)t+G(0)\leq  G(t)\leq C_7t^2+F(0)t+G(0),
\end{eqnarray}
where $C_6=\min\{4-n,n(\gamma-2)+2\}(E_k(0)+E_i(0))+(n-2)E_\delta(0)$ and $C_7=\max\{4-n,n(\gamma-2)+2\}(E_k(0)+E_i(0))+(n-2)E_\delta(0)$.
\end{Proposition}
\textbf{{\em Proof.}} In view of Lemma \ref{(100000)}, we have for (IEP)
\begin{eqnarray}\label{(36)}
     {\f {d^2}{dt^2}}G(t)\nonumber
     &=&IH(t)
     =2 E_k(t)+n(\gamma
-1)I(t)-{\f{\delta(n-2)}{2}} \int_{\mathbb{R}^n}
\rho \Phi dx\\
     &=&(4-n)E_k(t)+(n(\gamma-2)+2)I(t)+(n-2)C_0,
\end{eqnarray}
 and for (EP)
\begin{eqnarray}\label{(38)}
     {\f {d^2}{dt^2}}G(t)\nonumber
     &=&H(t)
     =2 E_k(t)+n(\gamma
-1)E_i(t)-{\f{\delta(n-2)}{2}}\int_{\mathbb{R}^n}
\rho \Phi dx\\
     &=&(4-n)E_k(t)+(n(\gamma-2)+2)E_i(t)+(n-2)E_\delta(0).
\end{eqnarray}

Therefore, by Corollary \ref{(100006)}, we can estimate ${\f {d^2}{dt^2}}G(t)$ as follows:\\
(I) For (IEP) with $\delta=-1$, we have
 \begin{equation}\label{(38.1)}
{\f {d^2}{dt^2}}G(t)\leq C_3, \ {\rm if}\  \gamma>2(1-{\f{1}{n }}),
 \end{equation}
and
\begin{equation}\label{(38.2)}
{\f {d^2}{dt^2}}G(t)\geq C_4.
 \end{equation}
(II) For (IEP) with $\delta=1$, we have
 \begin{equation}\label{(38.3)}
{\f {d^2}{dt^2}}G(t)\leq C_5.
 \end{equation}
(III) For (EP), we have
 \begin{equation}\label{(38.4)}
C_6\leq{\f {d^2}{dt^2}}G(t)\leq C_7.
 \end{equation}
 Integrating \eqref{(38.1)}-\eqref{(38.4)} over $[0,t]$, we get  \eqref{(33)}-\eqref{(37)}.

The following inequality is due to Chemin in \cite{Chemin}, one can also see \cite{JWX} for more details.

\begin{Lemma}\label{(100011)} For any $f\in L^1(\mathbb{R}^n,dx)\cap L^\gamma(\mathbb{R}^n,dx)\cap
L^1(\mathbb{R}^n,|x|^2dx)$, it holds that
\begin{eqnarray}\label{(39)}
   \parallel f \parallel_{L^1(\mathbb{R}^n,dx)}\leq C_8\parallel f \parallel_{L^\gamma(\mathbb{R}^n,dx)}^{\f{2\gamma}{(n+2)\gamma-n}}\parallel f \parallel_{L^1(\mathbb{R}^n,|x| ^2dx)}^{\f{n(\gamma-1)}{(n+2)\gamma-n}},\label{(25-)}
\end{eqnarray}
where $C_8=2|B_1|^{\f{2(\gamma-1)}{(n+2)\gamma-n}}=2({\f{\pi^{\f n 2}}{\Gamma ({\f n 2}+1)}})^{\f{2(\gamma-1)}{(n+2)\gamma-n}}$.
\end{Lemma}

Taking $f=\rho$ in Lemma \ref{(100011)}, we arrive at the lower bound
of $E_i(t) $ and $I(t)$, which read as:

\begin{Proposition}\label{(100012)}
 For (EP) and (IEP), we have

\begin{eqnarray}\label{(40)}
E_i(t)\geq {\f {C_{9}}{G(t)^{\f {n (\gamma-1)}{2}}}},
\end{eqnarray}
and

\begin{equation}\label{(41)}
  I(t)\geq {\f {C_{10}}{G(t)^{\f {n (\gamma-1)}{2}}}},
\end{equation}
respectively, where $C_{9}=({\f{\Gamma ({\f n 2}+1)}{(\pi)^{\f n
2}}})^{\gamma-1 } {\f{\exp ({\f {s_1}{c_\nu}})M(0)^{\f
{(n+2)\gamma-n}{2}}}{2^{\f{(n+2)\gamma-n}{2}}(\gamma-1)}}$ and
$C_{10}=({\f{\Gamma ({\f n 2}+1)}{(\pi)^{\f n 2}}})^{\gamma-1 }
{\f{M(0)^{\f
{(n+2)\gamma-n}{2}}}{2^{\f{(n+2)\gamma-n}{2}}(\gamma-1)}}$.
\end{Proposition}

The following are the crucial estimates in deriving the upper bound of internal energy.

\begin{Lemma}\label{(100013)} Assume $1<\gamma\leq1+{\f{2}{n }}$.\\
 (I) For (IEP) with $\delta=-1$,  it holds that
 \begin{eqnarray}\label{(42)}
{\f {d}{dt}}IJ_\delta(t)\leq
\left\{
\begin{array}{ll}
&{\f{2-n(\gamma-1)}{t+1}}IJ_\delta(t), \quad \quad\quad\quad\quad\quad \quad\quad\quad \quad\quad\quad\quad\quad n=3  \ {\rm or}\ 4,\\
& {\f{2-n(\gamma-1)}{t+1}}IJ_\delta(t)+{\f{C_1(n-4)}{2(t+1)^{{\f{2(n-2)}{n(\gamma-1)}}-1}}}IJ_\delta(t)^{\f{n-2}{n(\gamma-1)}}, n>4 \ {\rm and}\ \gamma>{\f{2n}{n+2}}.
\end{array}
 \right.
\end{eqnarray}
(II) For (IEP) with $\delta=1$, we have
\begin{eqnarray}\label{(43)}
{\f {d}{dt}}IJ_\delta(t)\leq
\left\{
\begin{array}{ll}
&{\f{5-3\gamma}{t+1}}IJ_\delta(t)+{\f{C_1}{2(t+1)^{{\f{5-3\gamma}{3(\gamma-1)}}-1}}}IJ_\delta(t)^{\f{1}{3(\gamma-1)}},  n=3\ {\rm and}\ \gamma>{\f{6}{5}} ,\\
& {\f{2-n(\gamma-1)}{t+1}}IJ_\delta(t), \quad\quad\quad\quad\quad\quad \quad\quad\quad\quad\quad n\geq4.
\end{array}
 \right.
\end{eqnarray}
(III) For (EP), the following estimate holds
\begin{eqnarray}\label{(44)}
{\f {d}{dt}}J_\delta(t)\leq{\f{2-n(\gamma-1)}{t+1}}J_\delta(t)+{\f{(n-4)\delta}{2}}(t+1) \int_{\mathbb{R}^n}
\rho \Phi dx-{\f{\delta}{2}}(t+1)^2{\f {d}{dt}}\int_{\mathbb{R}^n}
\rho \Phi dx.
\end{eqnarray}
\end{Lemma}
\textbf{{\em
Proof.}} In view of Lemma \ref{(100000)}, one can compute that for (IEP)
\begin{eqnarray}\label{(45)}
{\f {d}{dt}}IJ_\delta(t)=(2-n(\gamma-1))(t+1)I(t)+{\f{(n-4)\delta}{2}}(t+1) \int_{\mathbb{R}^n}
\rho \Phi dx
\end{eqnarray}
and for (EP)
\begin{eqnarray}\label{(46)}
{\f {d}{dt}}J_\delta(t)&=&(2-n(\gamma-1))(t+1)E_i(t)+{\f{(n-4)\delta}{2}}(t+1) \int_{\mathbb{R}^n}
\rho \Phi dx\nonumber\\
&&-{\f{\delta}{2}} (t+1)^2{\f {d}{dt}}\int_{\mathbb{R}^n}
\rho \Phi dx.
\end{eqnarray}
Due to Lemma \ref{(100001)}, if we regard
$$ G(t)-(t+1)F(t)+(t+1)^2 E_k(t)$$
as a quadratic function of $(t+1)$, since
$$\Delta=(F(t)^2-4G(t)E_k(t))\leq0,$$
 we have
$$G(t)-(t+1)F(t)+(t+1)^2 E_k(t)\geq0.$$
Consequently,
\begin{eqnarray}\label{(47)}
  E_i(t)\leq {\f{1}{(t+1)^2}}J(t),  \ \ I(t)\leq {\f{1}{(t+1)^2}}IJ_\delta(t).
\end{eqnarray}
Hence, by \eqref{(45)}, \eqref{(46)}, \eqref{(47)}  and Lemma \ref{(100003)}, we get the following estimates:\\
 (I) For (IEP) with $\delta=-1$, since $\int_{\mathbb{R}^n}\rho \Phi dx<0$, we have as $n=3$ or $n=4$
 \begin{eqnarray}
{\f {d}{dt}}IJ_\delta(t)\leq (2-n(\gamma-1))(t+1)I(t)\leq {\f{2-n(\gamma-1)}{t+1}}IJ_\delta(t)
\end{eqnarray}
 and when $n>4$ and $\gamma>{\f{2n}{n+2}}$
 \begin{eqnarray}\label{(47.5)}
{\f {d}{dt}}IJ_\delta(t)&=&(2-n(\gamma-1))(t+1)I(t)+{\f{(n-4)\delta}{2}}(t+1) \int_{\mathbb{R}^n}\rho \Phi dx\nonumber\\
&\leq&{\f{2-n(\gamma-1)}{t+1}}IJ_\delta(t)+{\f{C_1(n-4)}{2(t+1)^{{\f{2(n-2)}{n(\gamma-1)}}-1}}}IJ_\delta(t)^{\f{n-2}{n(\gamma-1)}}.
\end{eqnarray}
 Cases (II) and (III) are similar.

 It follows from  Lemma \ref{(100013)} that

\begin{Proposition}\label{(100014)}
Assume $1<\gamma\leq1+{\f{2}{n }}$.\\
(I) For (IEP) with $\delta=-1$. If $n=3$ or $n=4$, then
\begin{eqnarray}\label{(48)}
I(t)\leq{\f{C_{11}}{(t+1)^{n(\gamma-1)}}},
\end{eqnarray}
where $C_{11}=IJ(0)$.\\
(I) For (IEP) with $\delta=1$. If $n\geq4$, then
\begin{eqnarray}\label{(48)}
I(t)\leq{\f{C_{11}}{(t+1)^{n(\gamma-1)}}},
\end{eqnarray}
where $C_{11}=IJ(0)$.
\end{Proposition}

Now  we are ready to prove Theorem \ref{(10000)}--Theorem
\ref{(10002)}. We only give the proof of Theorem \ref{(10000)}, since the proof of Theorem \ref{(10001)} and Theorem \ref{(10002)} is similar.\\
\textbf{\textbf{Proof of Theorem \ref{(10000)}}.}
Suppose that the life span of the classical solution $T=+\infty$.  If one of (i) or (ii) occurs, we would deduce $G(t_0)<0$ for some time $t_0$, which contradicts the fact $G(t)\geq0$. We next  show the case (iii).

 By Proposition \ref{(100012)} and Proposition \ref{(100014)},  if
$2(1-{\f{1}{n }})<\gamma \leq1+{\f{2}{n }}$,  actually, n=3, then we have
\begin{equation}\label{(49)}
{\f {C_{10}}{G(t)^{\f {3 (\gamma-1)}{2}}}}\leq
I(t)\leq{\f{C_{11}}{(t+1)^{3(\gamma-1)}}},
\end{equation}
for all $t\geq0$.
In view of \eqref{(33)}, one has
\begin{equation}\label{(50)}
 G(t)\leq (3C_0+C_2)t^2+F(0)t+G(0),
\end{equation}
here we have used $C_3=(6-n)C_0+(4-n)C_2=3C_0+C_2$ since $2(1-{\f{1}{n }})<\gamma \leq1+{\f{2}{n }}$.
Substituting  \eqref{(50)} to  \eqref{(49)} yields
\begin{equation}\label{(51)}
{\f {C_{10}}{((3C_0+C_2)t^2+F(0)t+G(0))^{\f {3
(\gamma-1)}{2}}}}\leq{\f{C_{11}}{(t+1)^{3(\gamma-1)}}}.
\end{equation}
Let $t$ goes to infinity,  we get
\begin{equation}\label{(52)}
C_{10}\leq C_{11}(3C_0+C_2)^{\f {3
(\gamma-1)}{2}}.
\end{equation}
This leads to a contradiction with \eqref{(9)}. On the contrary,  if \eqref{(9)} holds,
then there  exists  a  time $T_1^{*}<\infty$, satisfying
\begin{equation}\label{(53)}
{\f {C_{10}}{((3C_0+C_2){T_1^{*}}^2+F(0){T_1^{*}}+G(0))^{\f {3
(\gamma-1)}{2}}}}>{\f{C_{11}}{(T_1^{*}+1)^{3(\gamma-1)}}},
\end{equation}
such that $T_1^*$ is the life span of the classical solution.
Indeed, one can solve out $T_1^{*}$ by  \eqref{(53)}. The proof of the
theorem  is finished.

\section*{Acknowledgement}
The author is deeply grateful to Professor Tao Luo  and Huihui Zeng for their invaluable suggestions, discussions and advices.

\end{document}